\documentclass[reqno, 12pt]{amsart}
\usepackage{amsmath,amsfonts,amsthm,amscd,upref,amstext}
\usepackage{amssymb}
\usepackage{subfigure}
\usepackage[all]{xy}
\usepackage{comment}
\usepackage{xcolor}
\usepackage{caption}
\usepackage{fullpage}
\usepackage[colorinlistoftodos]{todonotes}

\newcommand{\C}{\mathbb{C}}

\newcommand{\OO}{\mathcal{O}}

\newtheorem{prop}{Proposition}[section]
\newtheorem{thm}[prop]{Theorem}
\newtheorem{cor}[prop]{Corollary}
\newtheorem{ex}[prop]{Example}
\newtheorem{defn}[prop]{Definition}

\newtheorem{conjec}[prop]{Conjecture}
\newtheorem{rem}[prop]{Remark}
\newtheorem{lem}[prop]{Lemma}
\newtheorem{nota}[prop]{Notation}

\DeclareMathOperator{\edim}{edim}

\usepackage{tikz}

\usepackage{tikz-cd}
\usetikzlibrary{matrix}

\usepackage{lmodern} 
\usepackage{amsmath}
\usepackage{amssymb}
\usepackage{graphicx,scalerel}

\usepackage{scalerel,stackengine}
\stackMath

\newcommand\reallywidehat[1]{%
	\savestack{\tmpbox}{\stretchto{%
			\scaleto{%
				\scalerel*[\widthof{\ensuremath{#1}}]{\kern-.8pt\bigwedge\kern-.8pt}%
				{\rule[-\textheight/2]{1ex}{\textheight}}
			}{\textheight}%
		}{0.85ex}}%
	\stackon[1pt]{#1}{\tmpbox}%
}

\numberwithin{equation}{section}

\begin{document}
	
	\title{On the Gluing of germs of complex analytic spaces, Betti numbers and their structure}
	
	\date{}
	
	\author{T. H. Freitas}
	\address{Universidade Tecnol\'ogica Federal do Paran\'a, 85053-525, Guarapuava-PR, Brazil}
	\email{freitas.thf@gmail.com}
	
	\author{J. A. Lima}
	\address{Universidade Tecnol\'ogica Federal do Paran\'a, 85053-525, Guarapuava-PR, Brazil}
	\email{seyalbert@gmail.com}
	
	\thanks{The authors were partially supported  CNPq-Brazil Universal 421440/2016-3.}
	
	
	\begin{abstract}
		
		In this paper we introduce  new classes of gluing of  complex analytic spaces germs, called weakly large, large and strongly large. We  give a description of their Poincar\'e series and, as applications, we give numerical criteria to determine when these classes of gluing of germs of complex analytic spaces are smooth, singular, complete intersections and Gorenstein in terms of their Betti numbers. In particular, we show that the gluing of the same germ of complex analytic space along of any subspace is always a singular germ.  
		
	\end{abstract}
	
	\maketitle
	
	\section{Introduction}
	
	In modern algebraic geometry, gluing constructions are  a relevant topic of investigation by several authors over the years (for instance \cite{FREJOM}, \cite{paperformal}, \cite{gisri} and \cite{KS}). In the case that $(X,x)\subset (\C^n,x)$, $(Y,y)\subset (\C^m,y)$ and $(Z,z)\subset (\C^l,z)$ are germs of  complex analytic spaces, in \cite{FREJOM} 
	the  authors have shown that the gluing ${(X,x)} \sqcup_{(Z,z)} {(Y,y)}$ is also a germ of a complex analytic space, provided $\OO_{X,x}\to \OO_{Z,z}$ and $\OO_{Y,y}\to  \OO_{Z,z}$ are both surjective homomorphisms. Also, is given  the description of some algebraic/geometric and topological invariants such as the degree of a finite map germ, multiplicity and Milnor number.

	\begin{figure}[h!]
		\centering
		\includegraphics[width = 5.65cm, height = 4.5cm]{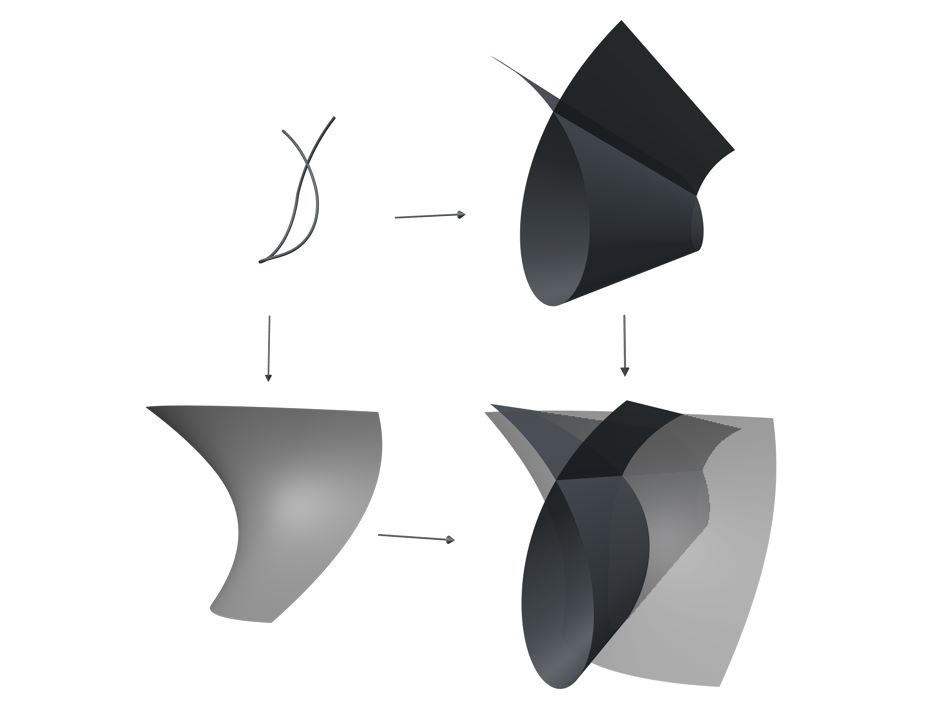}
		\caption*{Gluing of two surfaces along a curve }
		\label{figura1}
	\end{figure}
	
	The study of the structure of a germ of a complex analytic space deserves a special attention for the investigation  in the Singularity theory (\cite{gaffney}, \cite{giba}, \cite{MoVan}, \cite{nuno} and \cite{cidinha}). Some results concerning the Cohen-Macaulayness of the gluing ${(X,x)} \sqcup_{(Z,z)} {(Y,y)}$ are provided in \cite{FREJOM}, when $(Z,z)$ is a reduced point, and show that the structure of the gluing of germs of  analytic spaces  may have severe changes, depending on how this gluing is being made. For instance, the gluing of two Cohen-Macaulay surfaces can not be  Cohen-Macaulay, and the gluing of two germs of analytic spaces  that are complete intersection with isolated singularities (ICIS) is not always a complete intersection (see  \cite[Proposition 4.1, Theorem 4.3]{FREJOM}). When $(Z,z)$ is not a reduced point, results concerning the structure of the gluing (when is singular, smooth, complete intersection or Gorenstein) are not known.

	The main focus of the present paper is to define classes of gluing of germs of a complex analytic space, called weakly large, large and strongly large gluing,  and give numerical criteria to determine when it is smooth, singular, hypersurface, complete intersection and Gorenstein. The class of strongly large gluing contains, for instance,  the gluing  ${(X,x)} \sqcup_{(Z,z)} {(Y,y)}$ when $(Z,z)$ is a reduced point, and ${(X,x)} \sqcup_{(Z,z)} {(X,x)}$.  For this purpose, we give a description of their Poincar\'e series and the  Betti numbers in terms of the germs involved.

	The  crucial algebraic ingredient for this work is the notion of fiber product ring (pullback). This concept appears in several branches of mathematics and has been an ongoing topic of investigation in the category of rings (see \cite{AAM},   \cite{shirogoto}, \cite{paper} and \cite{nstv}).
	
	We briefly describe the contents of the paper. In Section 2,  we recall  the main definitions and results for the rest of the work. Section 3 is devoted to define the classes of weakly large, large and strongly large gluing of  complex analytic space germs, to give a shape of their Poincar\'e series (see Lemma \ref{lemgeral} and Theorem \ref{DKPoinFib}) and, as the main consequence, the description of their  Betti numbers (Corollaries \ref{cb1} and \ref{cormain2}).

	The last section contains the  applications of the paper. Actually, using the obtained Betti numbers, we give numerical criteria to determine when the classes of gluing defined in Section 2  are singular, hypersurfaces, complete intersections and Gorenstein. For instance, large gluing of germs of complex analytic spaces can be smooth, but any strongly large gluing is singular (Proposition \ref{rfp} and Theorem \ref{rfp1} (i)). 
	These results also illustrate that, despite the defined gluing classes has a subtle difference, for example,  the Betti numbers of strongly large gluing of complex analytic space germs provide a better understanding with respect to their structure.   
	As the main consequence of this section, we derive that the gluing ${(X,x)} \sqcup_{(Z,z)} {(X,x)}$ is always singular (see Corollary \ref{rfp01}).

	\section{Setup and Background}
	
	In this section we  recall the main concepts and results for the rest of the paper.  For the basic definitions  see \cite{shihoko} and \cite{Von}.
	
	\begin{defn}
		{\rm 
			Let $\Omega\subset \C^n$ be an open subset. A  closed subset  $X \subset \Omega$  is called {\it an analytic subset (or analytic set)} of $\Omega$ if for all $x\in X$, there is an open neighborhood $V\subset \Omega$ of $x$ and a finite set of analytic functions  $f_1,\dots,f_s\in \OO_n(\Omega)$ defined on $V$ such that
			$$X\cap V =\{x\in V \ \mid \  f_1(x)=\dots=f_s(x)=0\}.$$
		}
	\end{defn}
	
	\begin{defn}{\rm A {\it ringed space} $(X,\OO_X)$ is a Hausdorff topological space $X$ together with a sheaf of rings $\OO_X$. 
			In this case, $\OO_X$ is a sheaf of commutative rings on an  analytic set $X$. In order to simplify, we write $X$ for the pair $(X,\OO_X)$. In particular, if the stalk $\OO_{X,x}$ is a
			local ring for every $x\in X$, we call $(X,\OO_X)$ a {\it locally ringed space}.}
	\end{defn}

	A pair $(\varphi, \varphi^{\ast}): (X,\OO_X)\to (Y,\OO_Y)$ is called a {\it morphism of ringed spaces} if the map $\varphi: X\to Y$ is  continuous and $\varphi^*: \OO_Y\to \varphi_*\OO_X$ is a morphism of sheaves of rings. Also,
	$\varphi_*\OO_X$ is the sheaf of commutative rings given by $\varphi_*\OO_X(U)=\Gamma(\varphi^{-1}(U),\OO_X)$,
	for any open subset $U\subset Y$. 
	
	A morphism of locally ringed spaces $(X, \OO_X)$ and $(Y,\OO_Y)$ is a morphism of ringed spaces $(\varphi, \varphi^{\ast}): (X, \OO_X)\to (Y,\OO_Y)$ such that for all $x\in X$, the homomorphism
	$\varphi^*_x: \OO_{Y,\varphi(x)}\to \OO_{X,x}$ induced from $\varphi^*$ is a local homomorphism, i.e., ${\varphi^*}^{-1}(\mathfrak{m}_{X,x})=\mathfrak{m}_{Y,\varphi(x)}$.
	
	A morphism $(\varphi, \varphi^{\ast}): (X, \OO_X)\to (Y,\OO_Y)$ is an isomorphism if $\varphi$ is a homeomorphism and $\varphi^*$ is an isomorphism of sheaves of rings.
	
	From now on, the ringed space $(X,\OO_X)$ and the morphism of ringed spaces  $(\varphi,\varphi^{\ast})$ will be denoted by $X$ and  $\varphi$, respectively.
	
	\begin{lem}{\cite[Theorem 6.1.10]{Von}}\label{stalk}
		Let $\varphi: X \to Y$ be a morphism of $\C$-ringed spaces. Then $\varphi$ is an isomorphism if and only if $\varphi$ is a homeomorphism and $\varphi_x^*$ is an isomorphism for every $x\in X$.
	\end{lem}
	
	\begin{defn}
		{\rm 
			A ringed space 
			$(X,\OO_X)$ is called an {\it analytic space} if every $x\in X$ has a neighborhood $U$ such that $(U,\OO_X(U))$ is isomorphic to a local model $(V,\OO_V )$ as locally ringed spaces, i.e., $V$ is an analytic subset of an open set $\Omega\subset \C^n$ for some  $n$, and $\OO_V=(\OO_n({\Omega})/\mathcal{I}_V)|_V$.}
	\end{defn}
	
	\begin{defn}{\rm 
			On the set $\mathfrak{A}$ of pairs $(X,x)$ consisting of an analytic space $X$ and its point $x$, we define a relation $\sim$ as follows:\\
			$(X,x)\sim (Y,y) \Leftrightarrow$  there is a neighborhood $U\subset X$ of $x$, a neighborhood $V\subset Y$ of $y$ and an isomorphism $f:U\cong V$ such that
			$f(x)=y$.\\
			The relation becomes an equivalence relation; let the quotient set $\mathfrak{G}:=\mathfrak{A}/\sim$.
			An element of $\mathfrak{G}$ is called a {\it  germ of an analytic space}, denoted by $(X,x)$.
			
			A morphism of germs $(X,x)\to (Y,y)$ is a germ of an analytic spaces map $ X \to Y.$ For an open $U\subset X$, a point $x\in U$, and an analytic map $\varphi : U\to  Y$ with $\varphi(x)=y$, we denote the induced germ by $\varphi_x :(X,x)\to (Y,y)$.
		}
	\end{defn}
	
	\begin{rem}\label{quocient}{\rm 
			It should be noted that the elements of the stalks $\OO_{X,x}$ are seen as germs at $x$ of holomorphic functions on $X$. Each  germ is represented by a holomorphic function $f\in\OO_X(U)$, defined on an open neighborhood $U$ of $x$. Conversely, each $f\in \OO_X(U)$ defines a unique germ at $x\in U$, which is denoted by $f_x$. Hence, since
			$(X,\OO_X)$ is an analytic space and $x=(a_1,\dots,a_n)\in X\subset \Omega\subset \C^n$, one has the isomorphism
			$$\OO_{X,x}\cong \OO_{\C^n,x}/\mathcal{I}_{X,x}\cong \C\{x_1-a_1,\dots, x_n-a_n\}/\mathcal{I}_{X,x},$$
			\noindent where $\mathcal{I}_{X,x}=\{f_x\in \OO_{\C^n,x}\,\,|\,\,\exists \  f\in \OO_{\C^n}(U)\,\, {\rm representing }\,\, f_x \,\, {\rm and }\,\, f|_{U\cap X}=0 \}.$  Now, the fact that  $\OO_{\C^n,x}$ is Noetherian gives that the ideal $\mathcal{I}_{X,x}$ is finitely generated, and so there exists $f_1,\dots, f_k\in \OO_{\C^n,x}$  such that $\mathcal{I}_{X,x}=\langle f_1,\dots, f_k \rangle$. For this paper, $\mathcal{I}_{X,x}$ is an ideal that defines the germ $(X,x)$ of an  analytic space.  Note that $\OO_{X,x}$ is an analytic $\C$-algebra and is a local ring with maximal ideal $\mathfrak{m}_{X,x}=\{f\in \OO_{X,x} \ | \ f(x)=0\}$.
		}
	\end{rem}
	
	Set  $X\coprod Y$ as the co-product or disjoint union of sets $X$ and $Y$.
	
	\begin{defn}\label{copro}{\rm
			Let $\alpha: Z \to X$ and $\beta: Z \to Y$ be morphisms of ringed spaces. Set
			$$X \sqcup_Z Y= X\coprod Y/\sim, $$
			\noindent where the relation $\sim $ is generated by relations of the form $x\sim y$ ($x\in X$, $y\in Y$), provided there exists $z\in Z$ such that $\alpha(z) = x$ and $\beta(z)=y$.
			
			Namely, it is the smallest equivalence relation on $X\coprod Y$ such that after passing to the quotient $X\coprod Y/\sim$ the following square becomes commutative
			\begin{equation}\label{diagram0}
				\xymatrix{Z \ar[r]^{\alpha}  \ar[d]_{\beta} & X  \ar[d]^{f}\\
					Y  \ar[r]^{g}        & X\sqcup_ZY. }
			\end{equation}
			\noindent where  $f$ and $g$ are the continuous natural maps.}            
	\end{defn} 
	
	Since $(X,\OO_X)$, $(Y,\OO_Y)$ and $(Z,\OO_Z)$ are ringed spaces,  \cite[Proposition 2.2]{KS} provides that $(X\sqcup_ZY, \OO_{X\sqcup_ZY})$ is also a ringed space, and so $f$ and $g$ becomes morphisms of ringed spaces. Also, note that this  definition satisfies the universal property by \cite[Theorem 2.3]{KS}.
	
	Analogous to the morphisms of germs $\alpha_z: (Z,z) \to (X,\alpha(z))$ and $\beta_z: (Z,z) \to (Y,\beta(z))$, the previous definition  can be made for germs of analytic spaces $(X,x)$, $(Y,y)$ and $(Z,z)$, and denoted by $(X,\beta(z))\sqcup_{(Z,z)}(Y,\alpha(z))$. In the rest of the paper $(X,\beta(z))\sqcup_{(Z,z)}(Y,\alpha(z))$ will be denoted by $(X,x)\sqcup_{(Z,z)}(Y,y)$, where $\alpha(z)=x$ and $\beta(z)=y$. When a germ $(Z,z)$ is a reduced point, i.e, $(Z,z) =(z,z)$, we will denote $(Z,z)$ by $\{z\}$.
	
	Now, we recall an important definition for this paper.
	
	\begin{defn}\label{deffiber}{\rm The fiber product of homomorphisms  $\alpha_z^{\ast}:\OO_{X,x} \to \OO_{Z,z}$, $\beta_z^{\ast}:\OO_{Y,y} \to \OO_{Z,z}$ of $\C$-algebras   is defined by
			$$
			\OO_{X,x} \times_{\OO_{Z,z}} \OO_{Y,y}:=\{(s,t)\in  \OO_{X,x} \times \OO_{Y,y}  \ \mid \ \alpha_z^{\ast}(s)=\beta_z^{\ast}(t)  \}.$$ 
		}
	\end{defn}
	
	By \cite[Lemma 1.2]{AAM}, the fiber product is also a commutative and local ring  with  maximal ideal given by $\mathfrak{m} = \mathfrak{m}_{X,x} \times_{\mathfrak{m}_{Z,z}} \mathfrak{m}_{Y,y}$, where $\mathfrak{m}_{X,x}$,  $\mathfrak{m}_{Y,y}$ and $\mathfrak{m}_{Z,z}$ are the maximal ideals of $\OO_{X,x}$, $\OO_{Y,y}$ and $\OO_{Z,z}$, respectively. Also it is a subring of $\OO_{X,x} \times \OO_{Y,y}$ and   universal with respect to the  commutative diagram
	\begin{equation}\label{diafiber}
		\xymatrix{\OO_{X,x} \times_{\OO_{Z,z}} \OO_{Y,y} \ar[r]^{\pi_1}\ar[d]_{\pi_2} & \OO_{X,x}  \ar[d]^{\alpha_z^{\ast}}\\
			\OO_{Y,y}  \ar[r]^{\beta_z^{\ast}}        & \OO_{Z,z}, }
	\end{equation}
	where $\pi_1(s,t) = s$ and $\pi_2(s,t)=t$ are natural surjections. Also in \cite[Section 1 (1.0.3)]{AAM} and \cite[Lemma 2.1]{shirogoto} the authors  have shown that $\OO_{X,x} \times_{\OO_{Z,z}} \OO_{Y,y}$ is a  Noetherian local ring if both $\alpha_z^{\ast}$ and $\beta_z^{\ast}$ are surjective maps. It is important to realize that the assumptions over the maps are  crucial for the Noetherianess of the fiber product ring \cite[Example 2.9]{FREJOM}. Also,  if   $(X,x)$, $(Y,y)$ and $(Z,z)$ are germs of analytic spaces, then $\OO_{X,x} \times_{\OO_{Z,z}} \OO_{Y,y}$ is a reduced ring (\cite[Proposition 4.2.18]{Ela}). For the fiber product $\OO_{X,x} \times_{\OO_{Z,z}} \OO_{Y,y}$ we assume  that $\OO_{X,x}\neq \OO_{Z,z}\neq \OO_{Y,y}$. Note that every $\OO_{X,x}$-module (or $\OO_{Y,y}$-module) is an $\OO_{X,x} \times_{\OO_{Z,z}} \OO_{Y,y}$-module via Diagram \ref{diafiber}.

	Below we summarize some key results for this paper showed in \cite{FREJOM}. 
	
	\begin{lem}\label{Truke}\cite[Proposition 2.10]{FREJOM} Let  $\alpha: Z \to X$ and $\beta: Z \to Y$ be  holomorphic mappings of analytic spaces. Then,
		$${\OO}_{(X,\alpha(z))\sqcup_{(Z,z)}(Y,\beta(z))}\cong \OO_{X,\alpha(z)} \times_{\OO_{Z,z}} \OO_{Y,\beta(z)}.$$
	\end{lem}

	\begin{lem}\label{AlgebraAnalytic}\cite[Lemma 3.1 and Corollary 3.3(b)]{FREJOM}
		Let   $(X,x)\subset (\C^n,x)$, $(Y,y)\subset (\C^m,y)$ and $(Z,z)\subset (\C^l,z)$ be germs of analytic spaces such that $\OO_{X,x}\to \OO_{Z,z}$ and $\OO_{Y,y}\to  \OO_{Z,z}$ are both surjective homomorphisms.
		\begin{itemize}
			\item[(i)] Then, $\OO_{X,x} \times_{\OO_{Z,z}} \OO_{Y,y}$ is an analytic $\C$-algebra.
			
			\item[(ii)]  There is a germ $(W, w)$ and an isomorphism $\OO_{W,w}\to \OO_{X,x} \times_{\OO_{Z,z}} \OO_{Y,y}$ of local analytic $\C$-algebras. In particular $(X,x)\sqcup_{(Z,z)}(Y,y)\cong (W,w).$
			
		\end{itemize}
		
	\end{lem}

	\begin{thm}\label{maintoestruct}\cite[Theorem 3.4]{FREJOM}
		Let $X$, $Y$ and $Z$ be analytic spaces such that  $\OO_X\to \OO_Z$ and $\OO_Y\to  \OO_Z$ are both surjective homomorphisms.   Then,
		$X\sqcup_ZY$ is an analytic space.
	\end{thm}
	
	\section{Poincar\'e series and Betti numbers of gluing of Germs of analytic spaces}
	
	The main focus of this section is to define new classes of gluing of germs of complex analytic spaces, and give  the shape of their Poincar\'e series and Betti numbers. For this purpose, two important definitions are necessary:
	\begin{defn}{\rm Let $(X,x)\subset (\C^n,x)$  be a germ of an analytic space. Let $M$ be a finitely generated $\OO_{X,x}$-module. The Poincar\'e series of $M$ is given by
			$$P_{M}^{\OO_{X,x}}(t):= \sum_{i\geq 0}\dim_{\C} {\rm Tor}_i^{\OO_{X,x}}\left(M, \C \right) t^i,$$
			where $\C:= \frac{\OO_{X,x}}{\mathfrak{m}_{X,x}}$ is the residue field. The number $\beta_i^{\OO_{X,x}}(M):=\dim_{\C} {\rm Tor}_i^{\OO_{
					X,x}}\left(M, \C\right)$ is called $i$-th Betti number of $M$.   Let $I$ be an ideal of $\OO_{X,x}.$ The Poincar\'e series of $\OO_{X,x}/I$ is denoted  by
			$$P_{(Z,z)}^{(X,x)}(t):=P_{\OO_{X,x}/I}^{\OO_{X,x}}(t),$$
			where $(Z,z)$ is  subspace  of $(X,x)$ defined by  the reduced ideal $I$ of $\OO_{X,x}$. The {\it $i$-th Betti number} of $(Z,z)$ is defined by
			$\beta_i^{(X,x)}(Z,z):=\beta_i^{\OO_{X,x}}\left(\frac{\OO_{X,x}}{I}\right).$
		}
	\end{defn}
	
	\begin{rem}\label{rem3.3}{\rm
			Let $(Z,z)$ be a  subspace  of $(X,x)$ defined by  the reduced ideal $I$. Set $\mu(Z,z)$ as the minimal number of generators of $\OO_{X,x}/I$. Then 
			$$P^{(X,x)}_{(Z,z)}(t)=\mu(Z,z)+tP^{(X,x)}_{(\Omega_1, \omega_1)}(t),$$ where $(\Omega_1, \omega_1)$ is the subspace that represents the first syzygy  of   $\OO_{X,x}/I$  over $\OO_{X,x}$ (see \cite{DK75}).
		} 
	\end{rem}
	
	The next definition of motivated by the work of Levin  \cite{levin}.
	
	\begin{defn}\label{largemap}{\rm Let   $f: (Y,y) \to (X,x)$ be a morphism of germs of complex analytic spaces, such that the induced map $\OO_{X,x}\to \OO_{Y,y}$ is a surjective homomorphism.  Then $f$ is said to be {\it large}  provided, for any $(Z,z)$ subspace of $(Y,y)$ considered as a subspace of $(X,x)$, the following equality happens
			$$P_{(Z,z)}^{(X,x)}=P_{(Z,z)}^{(Y,y)}P_{(Y,y)}^{(X,x)}.$$ 
			
		}
	\end{defn}

	Now, we are able to define new classes of gluing of germs of analytic spaces.
	
	
	\begin{defn}\label{propestrelanotacao} {\rm  
			Let   $(X,x)\subset (\C^n,x)$, $(Y,y)\subset (\C^m,y)$ and $(Z,z)\subset (\C^l,z)$ be germs of analytic spaces such that $\OO_{X,x}\to \OO_{Z,z}$ and $\OO_{Y,y}\to  \OO_{Z,z}$ are both surjective homomorphisms.

			\begin{itemize} 
				\item[(i)] We say that the gluing $(X,x)\sqcup_{(Z,z)}(Y,y)$ is {\it weakly large}, provided   
				$$P_{(K,k)}^{(X,x)\sqcup_{(Z,z)}(Y,y)}(t)= P_{(K,k)}^{(X,x)}(t)P_{(X,x)}^{(X,x)\sqcup_{(Z,z)}(Y,y)}(t)\,\, ,  $$
				where $(K,k)$ is the subspace of $(X,x)$ that represents  the kernel of the map $\alpha_z^\ast$ as $\OO_{X,x}$-module (see Diagram \ref{diafiber}).  

				\item[(ii)] The gluing $(X,x)\sqcup_{(Z,z)}(Y,y)$ is called {\it large}  provided the map $f$ is large (see Diagram \ref{diagram0} and Definition  \ref{largemap}). In addition, if  the map $g$ is also large, we call the gluing $(X,x)\sqcup_{(Z,z)}(Y,y)$ as {\it strongly large} gluing of germs of analytic spaces.  
			\end{itemize}
		}
	\end{defn}
	
	It easy to see that  every strongly large gluing is large and therefore weakly large. The next example and remark show that these new classes of gluing of germs of analytic spaces are non-empty and contain interesting types of gluing.

	\begin{ex} Let   $(X,x)\subset (\C^n,x)$, $(Y,y)\subset (\C^m,y)$ and $(Z,z)\subset (\C^l,z)$ be germs of analytic spaces such that $\OO_{X,x}\to \OO_{Z,z}$ and $\OO_{Y,y}\to  \OO_{Z,z}$ are both surjective homomorphisms. 
		
		\begin{itemize}
			\item[(i)] If the germ $(Z,z)$ is a reduced point, the gluing $(X,x)\sqcup_{\{z\}}(Y,y)$ is strongly large. In fact, by \cite[Proposition 3.1]{L81} the  maps $f$ and $g$ are large. 
			
			\item[(ii)] Suppose that there are surjective ring homomorphism $\OO_{X,x}\to \OO_{Y,y}$ and the kernel of $\OO_{Y,y}\to  \OO_{Z,z}$ is a weak complete
			intersection ideal in $\OO_{X,x}$. Then the gluing  $(X,x)\sqcup_{(Z,z)}(Y,y)$ is  large  \cite[Theorem 3.12]{ramati}. 
		\end{itemize}
	\end{ex}
	
	\begin{rem}\label{truquemap}{\rm
			If we  assume  germs    $(X,x)\subset (\C^n,x)$, $(Y,y)\subset (\C^m,y)$ and $(Z,z)\subset (\C^l,z)$ of analytic spaces such that $\OO_{Y,y}\to \OO_{X,x}\to \OO_{Z,z}$ are both surjective homomorphisms,  \cite[3.11]{ramati} gives that the map $g$ (Diagram \ref{diagram0}) is large and therefore  the gluing $(X,x)\sqcup_{(Z,z)}(Y,y)$ is large. In particular one has that the gluing $(X,x)\sqcup_{(Z,z)}(X,x)$ is strongly large.
		}
	\end{rem}

	\begin{nota}{\rm  Throughout this paper, in order to use the structural results given in Lemma \ref{Truke}, Lemma \ref{AlgebraAnalytic} and Theorem \ref{maintoestruct},  we  assume  germs of analytic spaces   $(X,x)\subset (\C^n,x)$, $(Y,y)\subset (\C^m,y)$ and $(Z,z)\subset (\C^l,z)$  such that $\OO_{X,x}\to \OO_{Z,z}$ and $\OO_{Y,y}\to  \OO_{Z,z}$ are both surjective homomorphisms.} 
	\end{nota}
	
	We pose the following conjecture:
	
	\begin{conjec} {\rm Every gluing ${(X,x)\sqcup_{(Z,z)}(Y,y)}$ of complex analytic space germs  is large.}
	\end{conjec}
	
	The next result is a key ingredient for the rest of the paper and  shows the explicit shape of the Poincar\'e series of certain gluing of germs of analytic spaces. 
	
	\begin{lem}\label{lemgeral} Let   $(X,x)\subset (\C^n,x)$, $(Y,y)\subset (\C^m,y)$ and $(Z,z)\subset (\C^l,z)$ be germs of analytic spaces.   
		\begin{itemize}
			\item[(i)] If the gluing ${(X,x)\sqcup_{(Z,z)}(Y,y)}$ is weakly large, then
			$$P^{{(X,x)\sqcup_{(Z,z)}(Y,y)}}_{(X,x)}(t)=\frac{1-P_{{(Y,y)}}^{{(X,x)} \sqcup_{(Z,z)} {(Y,y)}}(t)}{1-P^{{(X,x)} }_{{(Z,z)}}(t)}.$$

			\item[(ii)] Suppose that the gluing ${(X,x)} \sqcup_{(Z,z)} {(Y,y)}$   is large. If $(W,w)$ is a subspace of $({X,x})$, then
			$$P^{{(X,x)} \sqcup_{(Z,z)} {(Y,y)}}_{(W,w)}(t)=\frac{P^{{(X,x)}}_{(W,w)}(t)\left(1- P^{{(X,x)} \sqcup_{(Z,z)} {Y,y}}_{{(Y,y)}}(t)\right)}
			{1-P^{{(X,x)}}_{{(Z,z)}}(t)}.$$
		\end{itemize}
	\end{lem}
	\begin{proof}
		(i) The exact sequence
		\begin{equation}\label{th31}
			0\longrightarrow {\rm ker}(\alpha_z^\ast)\longrightarrow \OO_{X,x}\times_{\OO_{Z,z}}\OO_{Y,y}\overset{\pi_2}{\longrightarrow} \OO_{Y,y}\longrightarrow 0
		\end{equation}
		and Remark \ref {rem3.3} gives
		\begin{equation}\label{th320}
			P_{(Y,y)}^{{(X,x)} \sqcup_{(Z,z)} {(Y,y)}}(t)=1+tP^{{(X,x)} \sqcup_{(Z,z)} {(Y,y)}}_{(K,k)}(t),
		\end{equation}
		where $(K,k)$ is the subspace of $(X,x)$ that represents  the kernel of the map $\alpha_z^\ast$.
		Since  ${(X,x)\sqcup_{(Z,z)}(Y,y)}$ is weakly large, one obtains
		\begin{equation}\label{th330}
			tP^{{(X,x)} \sqcup_{(Z,z)} {(Y,y)} }_{(K,k)}(t)= tP_{(K,k)}^{(X,x)}(t)P_{(X,x)}^{{(X,x)} \sqcup_{(Z,z)} {(Y,y)}}(t)=\left(P_{(Z,z)}^{(X,x)}(t)-1\right)P^{{(X,x)} \sqcup_{(Z,z)} {(Y,y)}}_{(X,x)}(t),
		\end{equation}
		where the last equality follows by the exact sequence
		\begin{equation}
			0\longrightarrow {\rm ker}(\alpha_z^\ast)\longrightarrow \OO_{X,x}\overset{{\rm ker}(\alpha_z^\ast)}{\longrightarrow} \OO_{Z,z}\longrightarrow 0
		\end{equation}
		and Remark \ref{rem3.3}.
		Hence (\ref{th320}) and (\ref{th330}) provide  \begin{equation}\label{th34}
			P_{(Y,y)}^{{(X,x)} \sqcup_{(Z,z)} {(Y,y)}}(t)=1+\left(P_{(Z,z)}^{(X,x)}(t)-1\right)P^{{(X,x)} \sqcup_{(Z,z)} {(Y,y)}}_{(X,x)}(t),
		\end{equation}
		and therefore
		\begin{equation}\label{23}
			P^{{(X,x)} \sqcup_{(Z,z)} {(Y,y)}}_{(X,x)}(t)=\frac{1-P_{(Y,y)}^{{(X,x)} \sqcup_{(Z,z)} {(Y,y)}}(t)}{1-P^{(X,x)}_{(Z,z)}(t)}.   
		\end{equation}
		
		(ii)  Since the gluing ${(X,x)} \sqcup_{(Z,z)} {(Y,y)}$   is large, by definition it is also weakly large. Hence,  multiplying both sides of (\ref{23})   by $P^{(X,x)}_{(W,w)}(t)$,  one has
		$$P_{(W,w)}^{{(X,x)} \sqcup_{(Z,z)} {(Y,y)}}(t)=P^{(X,x)}_{(W,w)}(t)+P^{(X,x)}_{(Z,z)}(t)P_{(W,w)}^{{(X,x)} \sqcup_{(Z,z)} {(Y,y)}}(t)-P^{{(X,x)} \sqcup_{(Z,z)} {(Y,y)}}_{(Y,y)}(t)P_{(W,w)}^{(X,x)}(t).$$
		Therefore
		$$P_{(W,w)}^{{(X,x)} \sqcup_{(Z,z)} {(Y,y)}}(t)=\frac{P_{(W,w)}^{(X,x)}(t)\left(1-P_{(Y,y)}^{{(X,x)} \sqcup_{(Z,z)} {(Y,y)}}(t)\right)}{1-P_{(Z,z)}^{(X,x)}(t)}.$$
		
	\end{proof}
	

	As a consequence, we derive a formula to compute the Betti numbers of any subspace of  the complex analytic germ $({X,x})$ as a subspace of the large gluing ${(X,x)} \sqcup_{(Z,z)} {(Y,y)}$.
	
	For the next two results, in order to simplify the notation, let $(\mathcal{V},v)$ denote the gluing of germs ${(X,x)} \sqcup_{(Z,z)} {(Y,y)}$ and $\beta_{i}^{T}(U):= \beta_{i}^{(T,t)}(U,u)$, for any germs  $(T,t)$ and $(U,u)$.
	
	\begin{prop}\label{bsp} Let $(X,x)\subset (\C^n,x)$, $(Y,y)\subset (\C^m,y)$ and $(Z,z)\subset (\C^l,z)$ be  germs of analytic spaces. Suppose that the gluing  $(\mathcal{V},v)$ is large. Then, for any $(W,w)$ subspace of $({X,x})$,
		$$\sum_{i=0}^{j-1}\beta_i^{\mathcal{V}}(W)\beta_{j-i}^{X}({Z})=\sum_{i=0}^{j-1}\beta_i^{X}(W)\beta_{j-1}^{\mathcal{V}}{(Y)},$$
		for each $j\geq 1$ positive integer.
	\end{prop}
	\begin{proof}
		Since the gluing $(\mathcal{V},v)$ is large,  Lemma \ref{lemgeral}  (ii) provides 
		$$P^{(\mathcal{V},v)}_{(W,w)}(t)=\frac{P^{{(X,x)}}_{(W,w)}(t)\left(1- P^{(\mathcal{V},v)}_{{(Y,y)}}(t)\right)}
		{1-P^{({X,x})}_{({Z,z})}(t)}.$$

		Set  $P_{(W,w)}^{({\mathcal{V},v})}(t)=\displaystyle\sum_i\beta_i^{{\mathcal{V}}}(W)t^i,$  $ 
		P_{(W,w)}^{(X,x)}(t)=\displaystyle\sum_i\beta_i^{X}(W)t^i,$  $
		P_{(Z,z)}^{({X,x})}(t)=\displaystyle\sum_i\beta_i^{X}(Z)t^i$   and
		
		$ P_{({Y,y})}^{(\mathcal{V},v)}(t)=\displaystyle\sum_i\beta_i^{{\mathcal{V}}}(Y)t^i.$ The previous equality yields
		\begin{equation}\label{ps1}
			\sum_i\beta_i^{{\mathcal{V}}}(W)t^i\left(1-\sum_i\beta_i^{{X}}({Z})t^i\right)=\sum_i\beta_i^{{X}}(W)t^i\left(1-\sum_i\beta_i^{{\mathcal{V}}}({Y})t^i\right).
		\end{equation}
		Note that
		$$
		\sum_i\beta_i^{{\mathcal{V}}}(W)t^i\left(1-\sum_i\beta_i^{{X}}(Z)t^i\right)=\sum_i\beta_i^{{\mathcal{V}}}(W)t^i- \sum_i\beta_i^{{\mathcal{V}}}(W)t^i\sum_i\beta_i^{X}(Z)t^i$$
		
		$$= \sum_i\beta_i^{{\mathcal{V}}}(W)t^i-\sum_{j\geq 0}\left(\sum_{i=0}^{j}\beta_i^{{\mathcal{V}}}(W)\beta_{j-i}^{X}(Z)\right)t^j.$$
		Similarly, the right side of equality  (\ref{ps1}) gives 
		{\footnotesize\begin{equation}\label{ps2}
				\sum_i\beta_i^{{\mathcal{V}}}(W)t^i-\sum_{j\geq 0}\left(\sum_{i=0}^{j}\beta_i^{{\mathcal{V}}}(W)\beta_{j-i}^{X}(Z)\right)t^j=\sum_i\beta_i^{{X}}(W)t^i-\sum_{j\geq 0}\left(\sum_{i=0}^{j}\beta_i^{{X}}(W)\beta_{j-i}^{{\mathcal{V}}}(Y)\right)t^j.
		\end{equation}}
		Therefore, for each $j\geq 1$,
		$$\beta_j^{{\mathcal{V}}}(W)-\sum_{i=0}^{j}\beta_i^{\mathcal{V}}(W)\beta_{j-i}^{{X}}(Z)=\beta_j^{X}(W)-\sum_{i=0}^{j}\beta_i^{X}(W)\beta_{j-i}^{{\mathcal{V}}}(Y).$$
		The fact $\beta_{0}^{X}(Z)=1=\beta_0^{{\mathcal{V }}}({Y})$ furnishes
		$$\sum_{i=0}^{j-1}\beta_i^{{\mathcal{V}}}(W)\beta_{j-i}^{X}(Z)=\sum_{i=0}^{j-1}\beta_i^{X}(W)\beta_{j-1}^{{\mathcal{V}}}(Y),$$ for each $j\geq 1$ positive integer and therefore, the desired conclusion follows.
	\end{proof}
	
	\begin{rem}\label{bettinaozero}{\rm  It is important to realize  that $\beta_1^{X}(Z)\neq 0\neq \beta_1^{Y}(Z)$, because otherwise, for instance, if $\beta_1^{X}(Z):=\beta_1^{\OO_{X,x}}(\OO_{Z,z})=0$, then $\OO_{Z,z}$ is a free $\OO_{X,x}$-module. The surjective map $\OO_{X,x}\stackrel{\alpha^{\star}_z}\to \OO_{Z,z}$  and the fact that $\OO_{Z,z}=\OO_{X,x}^{\oplus r}$, implies that $r=1$, (i.e., $\OO_{X,x}=\OO_{Z,z}$). This is a contradiction because  $\OO_{X,x}\neq \OO_{Z,z}$.}
	\end{rem}

	\begin{cor}\label{cb1}
		Let $(X,x)\subset (\C^n,x)$, $(Y,y)\subset (\C^m,y)$ and $(Z,z)\subset (\C^l,z)$ be  germs of analytic spaces. Suppose that the gluing  ${(X,x)} \sqcup_{(Z,z)} {(Y,y)}$ is large. Then, for any $(W,w)$ subspace of $({X,x})$,
		\begin{itemize}
			\item[(i)] $\beta_0^{\mathcal{V}}(W)=\displaystyle\frac{\beta_0^{X}(W)\beta_1^{\mathcal{V}}(Y)}{\beta_1^{X}(Z)}.$
			\item[(ii)]
			{\footnotesize{$\displaystyle\beta_1^{\mathcal{V}}(W)=\frac{1}{\beta_1^{X}(Z)}\left[\frac{\beta_1^{\mathcal{V}}(Y)\left(\beta_1^{X}(W)\beta_1^{X}(Z)-\beta_0^{X}(W)\beta_2^{X}(Z)\right)}{\beta_1^{X}(Z)}+\beta_0^{X}(W)\beta_2^{\mathcal{V}}(Y)\right].$}}
			
			\item[(iii)]{\tiny{$\beta_2^{\mathcal{V}}(W)=\displaystyle\frac{\beta_0^{X}(W)}{\beta_1^{X}(Z)}\left[\beta_3^{\mathcal{V}}(Y)+\frac{\beta_1^{\mathcal{V}}(Y)}{\beta_1^{X}(Z)}\Bigl(\beta_2^{X}(Z)k-\beta_3^{X}(Z)\Bigr)-k\beta_2^{\mathcal{V}}(Y)\right]+\frac{\beta_1^{X}(W)}{\beta_1^{X}(Z)}\Bigl[\beta_2^{\mathcal{V}}(Y)-k\beta_1^{\mathcal{V}}(Y)\Bigr]+\frac{\beta_2^{X}(W)\beta_1^{\mathcal{V}}(Y)}{\beta_1^{X}(Z)},$}}
			where $k=\displaystyle\frac{\beta_2^{X}(Z)}{\beta_1^{X}(Z)}.$
			
		\end{itemize}
	\end{cor}
	\begin{proof}
		(i)   By Proposition \ref{bsp} in the case $j=1$, one has  
		$$\beta_0^{\mathcal{V}}(W)\beta_1^{X}(Z)=\beta_0^{X}(W)\beta_1^{\mathcal{V}}(Y),$$
		which implies that $$\beta_0^{\mathcal{V}}(W)=\displaystyle\frac{{\beta_0^X(W)\beta_1^{\mathcal{V}}}(Y)}{\beta_1^{X}(Z)}.$$
		This gives (i). The proof of  (ii) and (iii) follows analogous by taking $j=2$ and $j=3$ in Proposition \ref{bsp}, respectively, together with the fact obtained in (i).
	\end{proof}
	
	\begin{thm} \label{DKPoinFib} Let  ${(X,x)} \sqcup_{(Z,z)} {(Y,y)}$ be the  gluing of the germs of analytic spaces $(X,x)\subset (\C^n,x)$, $(Y,y)\subset (\C^m,y)$ and $(Z,z)\subset (\C^l,z)$,
		satisfying one of the following conditions:
		\begin{itemize}
			\item[(i)] ${(X,x)} \sqcup_{(Z,z)} {(Y,y)}$ is weakly large and there is a surjective map $\OO_{Y,y} \twoheadrightarrow \OO_{X,x}$. 
			
			\item[(ii)] ${(X,x)} \sqcup_{(Z,z)} {(Y,y)}$ is strongly large.
		\end{itemize}
		If $(W,w)$ is a subspace of $(Y,y)$,   the Poincar\'e series of $(W,w)$ as a subspace of the gluing ${(X,x)} \sqcup_{(Z,z)} {(Y,y)}$ is given by
		$$P^{{(X,x)} \sqcup_{(Z,z)} {(Y,y)}}_{(W,w)}(t)=\frac{P^{{(Y,y)}}_{(W,w)}(t)P^{(X,x)}_{{(Z,z)}}(t)}
		{P^{{(X,x)}}_{{(Z,z)}}(t)+P^{{(Y,y)}}_{{(Z,z)}}(t)-P^{{(X,x)}}_{{(Z,z)}}(t)P^{{(Y,z)}}_{{(Z,z)}}(t)}.$$
	\end{thm} 
	\begin{proof}
		(i) Note that Lemma \ref{lemgeral} (i) (see (\ref{th34})) furnishes
		\begin{equation}\label{230}
			P_{(Y,y)}^{{(X,x)} \sqcup_{(Z,z)} {(Y,y)}}(t)=1+P_{({Z,z})}^{({X,x})}(t)P^{{(X,x)} \sqcup_{(Z,z)} {(Y,y)}}_{({X,x})}(t)-P^{{(X,x)} \sqcup_{(Z,z)} {(Y,y)}}_{(X,x)}(t).
		\end{equation}
		
		From the exact sequence 
		\begin{equation}
			0\longrightarrow {\rm ker}(\beta_z^\ast)\longrightarrow \OO_{X,x}\times_{\OO_{Z,z}}\OO_{Y,y}{\longrightarrow} \OO_{X,x}\longrightarrow 0,
		\end{equation}
		similarly  to the proof of Lemma \ref{lemgeral} (i),
		one obtains
		{\tiny\begin{equation}\label{th3300}
				P_{(X,x)}^{{(X,x)} \sqcup_{(Z,z)} {(Y,y)}}(t)=1+tP^{{(X,x)} \sqcup_{(Z,z)} {(Y,y)}}_{(K,k)}(t)=  1+ tP_{(K,k)}^{(Y,y)}(t)P_{(Y,y)}^{{(X,x)} \sqcup_{(Z,z)} {(Y,y)}}(t)= 1+(P_{({Z,z})}^{({Y,y})}(t)-1)P^{{(X,x)} \sqcup_{(Z,z)} {(Y,y)}}_{({Y,y})}(t),
		\end{equation}}
		
		\noindent where $(K,k)$ is the subspace of $(X,x)$ that represents  the kernel of the map $\alpha_z^\ast$, and the second equality follows by the hypothesis and  Remark \ref{truquemap}.

		Replacing   (\ref{th3300}) in   (\ref{230}) one has
		\begin{equation}\label{360}
			P_{{(Y,y)}}^{{(X,x)} \sqcup_{(Z,z)} {(Y,y)}}(t)=\frac{P_{(Z,z)}^{({X,x})}(t)}{P_{({Z,z})}^{({X,x})}(t)+ P_{({Z,z})}^{({Y,y})}(t)-P_{({Z,z})}^{({X,x})}(t)P_{({Z,z})}^{({Y,y})}(t)}.
		\end{equation}
		Again, by the hypothesis and 
		Remark \ref{truquemap},  multiplying both sides of equation (\ref{360}) by $P^{({Y,y})}_{(W,w)}(t)$ the desired conclusion follows.

		(ii) Since ${(X,x)} \sqcup_{(Z,z)} {(Y,y)}$ is strongly large, it is also weakly large. So, as in (\ref{230}),  \begin{equation}\label{2300}
			P_{(Y,y)}^{{(X,x)} \sqcup_{(Z,z)} {(Y,y)}}(t)=1+P_{(Z,z)}^{(X,x)}(t)P^{{(X,x)} \sqcup_{(Z,z)} {(Y,y)}}_{({X,x})}(t)-P^{{(X,x)} \sqcup_{(Z,z)} {(Y,y)}}_{(X,x)}(t).
		\end{equation}
		
		With an analogous argument used in (i) and a base change, it is possible to show  that 
		\begin{equation}\label{312}
			P_{(X,x)}^{{(X,x)} \sqcup_{(Z,z)} {(Y,y)}}(t)=1+P_{(Z,z)}^{(Y,y)}(t)P^{{(X,x)} \sqcup_{(Z,z)} {(Y,y)}}_{(Y,y)}(t)-P^{{(X,x)} \sqcup_{(Z,z)} {(Y,y)}}_{(Y,y)}(t),
		\end{equation}
		and therefore the statement is similarly obtained.
	\end{proof}
	
	\begin{cor}\label{amalgament} Let $(X,x)\subset (\C^n,x)$ and $(Z,z)\subset (\C^l,z)$ be  germs of analytic spaces. If $(W,w)$ is a subspace of $({X,x})$, then 
		$$P^{{(X,x)} \sqcup_{(Z,z)} {(X,x)}}_{(W,w)}(t)=\frac{P^{{(X,x)}}_{(W,w)}(t)}
		{2-P^{{(X,x)}}_{(Z,z)}(t)}.$$
	\end{cor}
	\begin{proof}
		The result is a consequence of Remark \ref{truquemap}  and Theorem \ref{DKPoinFib} (ii). 
	\end{proof}
	
	The next result shows the explicit shape of certain Betti numbers of the subspace $(W,w)$   of $(Y,y)$ seen as a  subspace 
	of the gluing ${(X,x)} \sqcup_{(Z,z)} {(Y,y)}$. We omit the proof because it is similar to  Corollary \ref{cb1}. 
	
	\begin{cor}\label{cormain2} Let  $(X,x)\subset (\C^n,x)$, $(Y,y)\subset (\C^m,y)$ and $(Z,z)\subset (\C^l,z)$ be germs of analytic spaces. Consider $(W,w)$  a subspace of $(Y,y)$. If the gluing ${(X,x)} \sqcup_{(Z,z)} {(Y,y)}$
		satisfies one of the conditions of Theorem \ref{DKPoinFib}, then
		\begin{itemize}
			\item[(i)] $\beta_0^{{(X,x)} \sqcup_{(Z,z)} {(Y,y)}}(W,w)=\beta_0^{(Y,y)}(W,w).$ 
			
			\item[(ii)] 
			$\beta_1^{{(X,x)} \sqcup_{(Z,z)} {(Y,y)}}(W,w)=\beta_0^{(Y,y)}(W,w)\beta_1^{(X,x)}(Z,z)+\beta_1^{(Y,y)}(W,w).$
			\item[(iii)] 
			$\beta_2^{{(X,x)} \sqcup_{(Z,z)} {(Y,y)}}(W,w)=\beta_0^{(Y,y)}(W,w)\beta_1^{(Y,y)}(Z,z)\beta_1^{(X,x)}(Z,z)+\beta_0^{(Y,y)}(W,w)\beta_2^{(X,x)}(Z,z)+\beta_1^{(Y,y)}(W,w)\beta_1^{(X,x)}(Z,z)+\beta_2^{(Y,y)}(W,w).$
		\end{itemize}
	\end{cor}

	\section{Applications}

	Again,  let $(\mathcal{V},v)$ denote the gluing ${(X,x)} \sqcup_{(Z,z)} {(Y,y)}$,  $\beta_{i}^{T}(U):= \beta_{i}^{(T,t)}(U,u)$, $\dim (T,t):= \dim (T)$ and the embedding dimension $\edim(T,t):= \edim (T)$, for any germs  $(T,t)$ and $(U,u)$. An important fact  for the rest of this section is that (\cite[Lemma 1.5 (1.5.2)]{AAM}) $$\dim (\mathcal{V})= {\rm max}\{\dim(X), \dim(Y)\}.$$

	As a consequence of the characterization of the Betti numbers of the gluing of germs of complex analytic spaces (Corollary \ref{cb1} (ii) and Corollary \ref{cormain2} (ii)), a formula for their embedding dimension is also provided.  
	\begin{cor}\label{edimfiber} Let  $(X,x)\subset (\C^n,x)$, $(Y,y)\subset (\C^m,y)$ and $(Z,z)\subset (\C^l,z)$ be germs of analytic spaces.
		
		\begin{itemize}
			\item[(i)] If the gluing $(\mathcal{V},v)$ is large, then
			$${\footnotesize{\displaystyle\edim({\mathcal{V}})=\frac{1}{\beta_1^{X}(Z)}\left[\frac{\beta_1^{\mathcal{V}}(Y)\left(\edim(X)\beta_1^{X}(Z)-\beta_2^{X}(Z)\right)}{\beta_1^{X}(Z)}+\beta_2^{\mathcal{V}}(Y)\right].}}$$
			
			\item[(ii)] If the gluing $(\mathcal{V},v)$
			satisfies one of the conditions of Theorem \ref{DKPoinFib}, then
			$$\edim({\mathcal{V}})=\beta_1^{X
			}(Z)+\edim(Y).$$
			
		\end{itemize}
	\end{cor}

	\begin{prop}\label{rfp} Let  $(X,x)\subset (\C^n,x)$, $(Y,y)\subset (\C^m,y)$ and $(Z,z)\subset (\C^l,z)$ be germs of analytic spaces such that $\dim (\mathcal{V})= \dim {(X)}.$ 
		Suppose that the gluing $(\mathcal{V},v)$
		is large.  Then $(\mathcal{V},v)$ is smooth if and only if the following equality holds
		$$\beta_1^{X}({Z})\left(\beta_1^{\mathcal{V}}{(Y)}\edim {(X)}-\dim {(X)}\beta_1^{{X}}({Z})+\beta_2^{\mathcal{V}}(T)\right)=\beta_1^{\mathcal{V}}(Y)\beta_2^{{X}}({Z}).$$
		
	\end{prop}
	\begin{proof}
		By definition, $(\mathcal{V},v)$ is smooth  if and only if $\edim (\mathcal{V})=\dim (\mathcal{V}).$ Hence, using Corollary \ref{cb1}  (ii), one has that $(\mathcal{V},v)$ is smooth if and only if
		\begin{equation}\label{ep312}
			\dim (X)=\frac{\beta_1^{\mathcal{V}}(Y)\left(\edim (X)\beta_1^{X}(Z)-\beta_2^{X}(Z)\right)}{\beta_1^{X}(Z)^2}+\frac{\beta_2^{\mathcal{V}}(Y)\beta_1^{X}(Z)}{\beta_1^{X}(Z)^2}.
		\end{equation}
		Solving  (\ref{ep312}) for $\beta_1^{X}(Z)$, one obtains that $(\mathcal{V},v)$  is smooth if and only if
		$$\beta_1^{X}(Z)\left(\beta_1^{\mathcal{V}}(Y)\edim (X)-\dim (X)\beta_1^{X}(Z)+\beta_2^{\mathcal{V}}(Y)\right)=\beta_1^{\mathcal{V}}(Y)\beta_2^{X}(Z),$$ and this shows the statement.
	\end{proof}
	
	As an immediate consequence of  Proposition \ref{rfp}, we derive the following.

	\begin{cor}
		Let  $(X,x)\subset (\C^n,x)$, $(Y,y)\subset (\C^m,y)$ and $(Z,z)\subset (\C^l,z)$ be germs of analytic spaces such that $\dim (\mathcal{V})= \dim (X).$ 
		Suppose that the gluing $(\mathcal{V},v)$
		is large and $\beta_1^{\mathcal{V}}(Y)=\beta_1^{{X}}({Y})=1$. Then $(\mathcal{V},v) $  is smooth if and only if $$\edim {(X)}-\dim ({X})=\beta_2^{{X}}({Z})-\beta_2^{\mathcal{V}}({Y}).$$
	\end{cor}

	\begin{prop}\label{p4.4} Let  $(X,x)\subset (\C^n,x)$, $(Y,y)\subset (\C^m,y)$ and $(Z,z)\subset (\C^l,z)$ be germs of analytic spaces. Suppose that $(\mathcal{V},v)$
		is large, $\beta_2^X(Z)=0$ and  $\beta_2^{\mathcal{V}}(Y)=\beta_1^{X}(Z).$ Then $(\mathcal{V},v)$ is a complete intersection if and only if 
		$$\edim ({X})+\dim ({X})= \frac{1}{m}\Bigl[\edim ({X})\beta_1^{\mathcal{V}}(Y)l-\beta_3^{\mathcal{V}}(Y)-\beta_2^{{X}}(k)\beta_1^{\mathcal{V}}{(Y)}\Bigr],$$
		where $l=\displaystyle\frac{3-m}{2m}$ with $m=\beta_1^{{X}}({Z}).$
		
	\end{prop}
	\begin{proof}
		By \cite[Theorem 7.3.3]{Avramov} (or \cite[Proposition 2.8.4 (3)]{HA}), $(\mathcal{V},v)$  is a complete intersection if and only if 
		\begin{equation}\label{IC}
			\beta_2^{\mathcal{V}}(0)=\binom{\beta_1^{\mathcal{V}}(0)}{2}+\beta_1^{\mathcal{V}}(0)-\dim (\mathcal{V}).
		\end{equation}
		The  equalities given  in Corollary \ref{cb1}  (ii)-(iii), the hypothesis $\beta_2^{\mathcal{V}}(Y)=\beta_1^{X}(Z)$, and the equation (\ref{IC}) yield the statement.
	\end{proof}
	
	It should be noted that the last results show that is difficult to have large gluing of germs of analytic  spaces  being smooth and complete intersection. For the  cases of Theorem \ref{DKPoinFib}, the next result  yields  a better understanding of their structure.
	
	\begin{thm}\label{rfp1} Let  $(X,x)\subset (\C^n,x)$, $(Y,y)\subset (\C^m,y)$ and $(Z,z)\subset (\C^l,z)$ be germs of analytic spaces such that $\dim (\mathcal{V})= \dim (Y).$  Suppose that  the gluing $(\mathcal{V},v)$ satisfies one of the conditions of Theorem \ref{DKPoinFib}. 
		\begin{itemize}
			\item[(i)] Then $(\mathcal{V},v)$  is singular.
			
			\item[(ii)] Suppose that $(\mathcal{V},v)$ is Cohen-Macaulay. Then $(\mathcal{V},v)$ is a hypersurface if and only if $(Y,y)$ is smooth and $\beta_1^{X}({Z})=1$.
			\item[(iii)] Suppose that $(X,x)$ is a complete intersection. Then $(\mathcal{V},v)$ is a complete intersection if and only if
			$$\frac{\beta_1^X(Z)^2+\beta_1^X(Z)}{\beta_1^X(Z)\beta_1^Y(Z)+\beta_2^X(Z)}=2.$$
			\item[(iv)] Suppose that $(\mathcal{V},v)$ is Cohen-Macaulay. If $\beta_1^{\mathcal{V}}(W)\leq \beta_0^{\mathcal{V}}(W)$ for some $(W, w)$ subspace of $(Y, y)$, then $(\mathcal{V}, v)$ is Gorenstein if and only if $(Y,y)$ is smooth.
		\end{itemize}
	\end{thm}
	\begin{proof}
		
		(i)  Suppose that $(\mathcal{V},v)$ is singular. Then, Corollary \ref{edimfiber} (ii) gives 
		$$\dim (Y)=\dim (\mathcal{V})= \edim(\mathcal{V})=\beta_1^{X}(Z)+\edim(Y).$$ 
		So $\beta_1^{X}(Z)=0$ because  $\edim(Y)\geq \dim (Y)$, which is a contradiction (Remark \ref{bettinaozero}).
		
		(ii) By (i), since $(\mathcal{V},v)$ is singular one has   $\edim(\mathcal{V})- \dim (\mathcal{V})>0$.  Hence $(\mathcal{V},v)$ is a hypersurface if and only if  $\edim(\mathcal{V})-{\rm depth} (\mathcal{V})=1$. Since  $(\mathcal{V},v)$ is Cohen-Macaulay by hypothesis,
		Remark \ref{edimfiber} furnishes
		$$\beta_1^{X}(Z)+\edim(Y)-{\rm dim} (Y)=1.$$
		The facts   $\beta_1^{X}(Z)\neq 0$ (Remark \ref{bettinaozero}) and $\edim(X)\geq \dim(X)$ yield that $(\mathcal{V},v)$ is a hypersurface if and only if   $\beta_1^{X}(Z)=1$ and $\edim(Y)={\rm dim}(Y)$ (i.e., $Y$ is smooth). 
		
		(iii) Set $d:=\dim (\mathcal{V})$. By \cite[Proposition 2.8.4 (3)]{HA}), $(\mathcal{V},v)$ is a  complete intersection if and only if
		\begin{equation}
			\label{equal} \beta_2^{\mathcal{V}}(0)= \binom{\overline{e}}{2}+\overline{e}-d,
		\end{equation}
		where  $\overline{e}:= \edim(\mathcal{V})=\beta_1^{X}({Z})+\edim({Y})$ (Corollary \ref{edimfiber} (ii)). 
		By Corollary \ref{cormain2} (iii) and the fact that $(X,x)$ is a complete intersection (\cite[Proposition 2.8.4 (3)]{HA}) yield   \begin{equation}\label{3}\beta_2^{\mathcal{V}}(0)=  \beta_1^{Y}(Z)\beta_1^{{X}}({Z})+\beta_2^{X}(Z)+e_2\beta_1^{X}(Z)+\binom{e_2}{2}+e_2-d,\end{equation} 
		where $e_2:=\edim(Y)$. Therefore, comparing  (\ref{equal}) and ($\ref{3}$) we obtains that $(\mathcal{V},v)$ is a complete intersection if and only if
		$$\beta_1^{X}(Z)^2+\beta_1^{X}(Z)=2\left(\beta_1^{Y}(Z)\beta_1^{X}(Z)+\beta_2^{X}(Z)\right).$$
		
		The desired conclusion follows, because  $\beta_1^{Y}(Z)\neq 0\neq \beta_1^{X}(Z)$ (Remark \ref{bettinaozero}).
		
		(iv) Suppose that $(Y,y)$ is smooth. Since   $\beta_1^{X}(Z)\neq 0$ (Remark \ref{bettinaozero}),  the hypothesis and Corollary \ref{cormain2} (i)-(ii) provide $\beta_1^{X}(Z)=1$ and $\beta_1^{Y}(W)=0$.  Therefore $(\mathcal{V},v)$ is Gorenstein by (ii). The converse immediately follows from \cite[Proposition 4.19]{shirogoto}.
	\end{proof}
	
	
	As mentioned in  Remark \ref{truquemap}, the gluing ${(X,x)} \sqcup_{(Z,z)} {(X,x)}$ is always strongly large. Since the dimension of $\dim{(X,x)} \sqcup_{(Z,z)} {(X,x)}$ and $\dim(X)$ are equal, as a consequence of Theorem \ref{DKPoinFib} we derive the following result.
	
	\begin{cor}\label{rfp01} Let  $(X,x)\subset (\C^n,x)$, and $(Z,z)\subset (\C^l,z)$ be germs of analytic spaces.   
		\begin{itemize}
			\item[(i)] Then ${(X,x)} \sqcup_{(Z,z)} {(X,x)}$ is singular.
			
			\item[(ii)] If ${(X,x)} \sqcup_{(Z,z)} {(X,x)}$ is Cohen-Macaulay, then ${(X,x)} \sqcup_{(Z,z)} {(X,x)}$ is a hypersurface if and only if $X$ is smooth and $\beta_1^X(Z)=1.$
			
			\item[(iii)] If $(X,x)$ is a complete intersection, then ${(X,x)} \sqcup_{(Z,z)} {(X,x)}$ is a complete intersection if and only if  
			$\beta_1^X(Z)=1$ and $\beta_2^X(Z)=0$. 
			
			\item[(iv)] Suppose that ${(X,x)} \sqcup_{(Z,z)} {(X,x)}$ is Cohen-Macaulay. If $\beta_1^{\mathcal{V}}(W)\leq \beta_0^{\mathcal{V}}(W)$ for some $(W, w)$ subspace of $(X, x)$, then $(\mathcal{V}, v)$ is Gorenstein if and only if $(X,x)$ is smooth.
			
		\end{itemize}
	\end{cor}
	\begin{proof} The proof of (i), (ii) and (iv) are immediate consequences of Theorem \ref{DKPoinFib} (i)-(ii)-(iv). For (iii), Theorem \ref{DKPoinFib} (iii) furnishes
		\begin{equation} \label{eqint0}\beta_1^X(Z)-\beta_1^X(Z)^2=2\beta_2^X(Z).\end{equation}
		Note that, if $\beta_1^X(Z)>1$, then left side of (\ref{eqint0}) is a negative number. Since $\beta_2^X(Z)\geq 0$, the equality (\ref{eqint0}) occurs if and only if $\beta_1^X(Z)=1$ and $\beta_2^X(Z)=0$ or $\beta_1^X(Z)=0$ and $\beta_2^X(Z)=0$. But $\beta_1^X(Z)\neq 0$ (Remark \ref{bettinaozero}), and therefore the result follows.      
	\end{proof}

	\noindent {\it Acknowledgements:}  
	The authors would like to thank Victor Hugo Jorge P\'erez and   Aldicio Jos\'e Miranda for the kind comments and suggestions for the improvement of the paper.

\end{document}